\documentclass{article}
\usepackage{amsmath,amssymb,amsthm,amsfonts,amscd,euscript,verbatim}

\usepackage{hyperref}

\usepackage{t1enc}
\theoremstyle{plain}
\newtheorem{theo}{Theorem}[section]

\newtheorem{pr}[theo]{Proposition}

\theoremstyle{remark}
\newtheorem{rema}[theo]{Remark}

\theoremstyle{definition}
\newtheorem{defi}[theo]{Definition}

\newcommand\smc{{SmCor}}

\newcommand\dmge{DM^{eff}_{gm}}
\newcommand\dmgm{DM_{gm}}
\newcommand\dme{DM_-^{eff}}

\newcommand\mg{\mathcal{M}}
\newcommand\obj{Obj}

\newcommand\id{id}
\newcommand\cu{\underline{C}}
\newcommand\du{\underline{D}}

\newcommand\z{\mathbb{Z}}
\newcommand\q{\mathbb{Q}}

\newcommand\p{\mathbb{P}}

\newcommand\ns{\{0\}}

\newcommand\chow{Chow}
\newcommand\cho{{Corr_{rat}}}
\newcommand\chowe{Chow^{eff}}
\newcommand\ab{Ab}
\newcommand\var{Var}
\newcommand\sv{SmVar}
\newcommand\spv{SmPrVar}

\newcommand\spe{\operatorname{Spec}\,}

\DeclareMathOperator\imm{\operatorname{Im}}
\DeclareMathOperator\co{\operatorname{Cone}}

\newcommand\znz{\z/n\z}
\newcommand\zop{{\mathbb{Z}[\frac{1}{p}]}}
\newcommand\zll{{\mathbb{Z}_{(l)}}}
\newcommand\ql{{\mathbb{Q}_{l}}}

\newcommand\gd{\mathfrak{D}}
\newcommand\gdp{\mathfrak{D}'}
\newcommand\gds{\mathfrak{D}_s}

\newcommand\eu{\underline{E}}
\newcommand\au{\underline{A}}

\newcommand\kar{\operatorname{Kar}}
\newcommand\cha{\operatorname{char}}

\newcommand\hrt{{\underline{Ht}}}
\newcommand\hw{{\underline{Hw}}}
\newcommand\w{{\mathfrak{w}}}
\newcommand\gdbr{\mathfrak{D}_{bir}}

\newcommand\ff{\mathbb{F}}
\newcommand\fq{\mathbb{F}_q}

\newcommand\dbm{D^b_m}
\newcommand\dsh{DSH}

\newcommand\dmcs{DM^c(S)}

\newcommand\wchow{{w_{Chow}}}

\DeclareMathOperator\adfu{\operatorname{AddFun}}

\begin{document}

 \title{Weight structures and motives; comotives,
 coniveau and Chow-weight spectral sequences, and mixed complexes of sheaves: a survey}

\author{M.V. Bondarko, St. Petersburg State University
\thanks{ 
The author gratefully acknowledges the support from Deligne 2004
Balzan prize in mathematics. The work is also supported by RFBR
(grants no. 08-01-00777a and 10-01-00287). } }

\maketitle

\tableofcontents

\section{Introduction}

This article is a  survey of author's results on Voevodsky's motives
and weight structures; yet it is supplied with detailed references. Weight structures are natural
counterparts of $t$-structures (for triangulated categories) introduced
 by the author in \cite{bws} (and also independently by D. Pauksztello in
 \cite{konk}). They allow to to construct {\it weight complexes},
{\it weight  filtrations}, and {\it weight spectral sequences}.
Partial cases of the latter are:  'classical' weight spectral
 sequences (for singular and \'etale cohomology), coniveau spectral
 sequences, and Atiyah-Hirzebruch spectral sequences
 (we mention all of these below). The details, proofs,
 and
several more results could be found in \cite{bws},
 \cite{bger}, and \cite{mymot} (we also mention certain results of \cite{bzp}, \cite{hebpo}, and \cite{brelmot}). We describe more motivation for the
 theory of weight structures, and define weight structures in
 \S\ref{sbws}.

 Though our 'main' weight structures will be defined on certain
 'motivic' categories, the author tried to make this survey accessible
 to readers that are rather interested in general triangulated categories
  (or possibly, the stable homotopy category in topology).
  Those readers may
freely ignore all definitions and results that are related with
algebraic geometry (and motives). On the other hand, the main
motivic results (see \S\ref{smmot}) could be understood without
knowing anything about weight structures (after \S\ref{smmot} a
'motivic' reader may proceed directly to \S\ref{swmot} to  find some
more motivation to study weight structures). Alternatively, it is
quite possible for any reader to read  section \S\ref{smmot} only
after studying the general theory of weight structures
(\S\S\ref{sbws}--\ref{sore}).

 The author chose not to pay much attention to the differential graded
 approach to motives in this  text; yet it is described in detail in
 \cite{mymot} and in \S6 of \cite{bws} (see also \cite{bev} and \S\ref{swc} below).

This text is based on the talks presented by the author at the
conferences
"Finiteness for motives and motivic cohomology"
 (Regensburg, 9--13th of February, 2009) and "Motivic homotopy theory"
 (M\"unster, 27-31st of July, 2009); yet some more recent topics are added.  The author is
deeply grateful to prof. Uwe Jannsen, prof. Eric Friedlander,
and to other organizers of
 these conferences for their efforts.

\section{Categoric notation; definitions of Voevodsky}\label{sdmot}

For a category $C,\ A,B\in\obj C$, we denote by
$C(A,B)$ the set of  $C$-morphisms from  $A$ to $B$.

Below
$B$ will be some additive category;
   $K^b(B)\subset K(B)$ will denote the homotopy
category of (bounded) $B$-complexes.

$\cu$ and $\du$ will be  triangulated categories;
for $f\in\cu (X,Y)$, $X,Y\in\obj\cu$, we will denote the third vertex
of (any) distinguished triangle $X\stackrel{f}{\to}Y\to Z$ by $\co(f)$.

For $D,E\subset \obj \cu$ we will write $D\perp E$ if $\cu(X, Y)=\ns$
 for all $X\in D,\ Y\in E$.

$\au$ will be an abelian category, $D(\au)$ is its derived category;
$H:\cu\to \au$ will usually be a cohomological functor (i.e.  it is
contravariant, and converts distinguished triangles into long exact
sequences in $\au$).

$\kar(B)$ for any $B$ will denote the Karoubization of $B$ i.e.
the category of 'formal images' of idempotents in $B$
(so $B$ is embedded into an idempotent complete category).

A full subcategory $C\subset B$ is called {\it Karoubi-closed} in $B$
if $C$ contains all
$B$-retracts of its objects; $\kar_{B}C$ will denote the smallest
Karoubi-closed subcategory of $B$ that contains $C$ (i.e. its objects
are all retracts of objects of $C$ that belong to $B$).

$\ab$ is the category of abelian groups.


Now we introduce our 'motivic' definitions; they could be especially
interesting to readers that are aware of 'classical' motives but do
not know much about Voevodsky's ones.

$k$ is our perfect base field. From time to time we will  have to assume that
either  $\operatorname{char} k=0$ or that we consider (co)motives
and cohomology with rational (or $\zop$-) coefficients.

$\spv\subset \sv\subset \var$ are the sets of (smooth projective)
 varieties over $k$.

\vspace{\baselineskip}

The definition of Voevodsky's motives starts from smooth
correspondences (see \cite{1}):
$\obj\smc =\sv$; $\smc(X,Y)=\z^{\{U\}}$: $U\subset
X\times Y$ is closed reduced, finite dominant over a
component of $X$. Compositions of morphisms are given
by a natural algebraic analogue of the composition of multi-valued functions.

\begin{rema}\label{rcoeff}
1. So, in contrast to the 'classical' definition, we consider
only those primitive correspondences (i.e. closed subvarieties
 of $X\times Y$ of a certain dimension) that are
finite over $X$. Note here that any 'classical' correspondence is rationally equivalent to some finite one.
 The advantage of finite correspondences is that the composition
is well-defined without factorizing modulo an equivalence relation.
This is very important! 

2. For any commutative associative ring
with a unit $R$ instead of $\smc$ one can consider a certain category $\smc_R$; in order to define it one should just
replace $\z^{\{U\}}$ by $R^{\{U\}}$ in the definition of $\smc(X,Y)$. This allows to construct a reasonable theory  of Voevodsky's motives with $R$-coefficients; see  \cite{vbook}. 

Usually one takes $R=\z$ or $R=\q$. In \cite{bzp} the author also considers certain intermediate coefficients rings. The case $R=\z/(n)$ (for $n>1$) is also interesting.
\end{rema}

Cartesian product of varieties yields tensor structure for $\smc$
(as well as for $K^b(\smc)$).

One can define (homological) Chow motives in terms of $\smc$.
One starts from the category of rational correspondences:
$\obj\cho=\spv$;
$\cho(X,Y)=\smc(X,Y)$/rational equivalence.

Now, one has $\chowe=\kar(\cho)$ (this yields a category that is
isomorphic to the 'classical' effective Chow motives).
 Formal tensor inversion of $\z(1)[2]$ (the Lefschetz motif i.e. the
 'complement' of a point to the projective line) yields the whole
 category $\chow$.

$\dmge$ is defined as the Karoubization of a certain localization of
$K^b(\smc)$ (so it is triangulated). Tensor inversion of
  $\z(1)[2]$ in it yields
 $\dmgm$.

We denote by $\mg$ the composition $\sv\to \smc\to K^b(\smc)\to \dmge$;
this defines  motives of smooth varieties. If $\operatorname{char} k=0$,
 in $\dmge$ there also exist motives and certain {\it motives with compact
 support} for arbitrary varieties.

Voevodsky constructed the following diagram of functors:

\begin{equation}\label{evo}
\begin{CD}
\chowe@>{}>>\chow\\
@VV{}V@VV{}V \\
\dmge@>{}>>\dmgm\end{CD}
\end{equation}

Here all arrows are full embeddings of additive categories.

In \S3.1 of \cite{1} Voevodsky also defined a certain triangulated
category $\dme\supset \dmge$.


\section{Main motivic results}
\label{smmot}

We list our main results. Assertions 1--6 require  $\cha k=0$
(yet see  Remark \ref{rmmot}(5) below). 


\begin{theo}\label{main}
1. In \S3 of \cite{mymot} $\dmge$ was described
 'explicitly' in terms of {\it twisted complexes} over a certain
  differential graded category
 $J$ (see \S2.4 of ibid.); the  objects of $J$ are
cubical Suslin complexes of smooth projective varieties.

2. This description is somewhat similar to (yet 'more convenient' than)
those
of Hanamura's motives (see \cite{h}). This allowed to compare
Voevodsky's motives with
Hanamura's ones: in \S4 of \cite{mymot} it was
proved that $\dmgm\q$ is anti-isomorphic to Hanamura's
motives.

3.  'Killing
all arrows of negative degrees' in the 'description' of $\dmge$
 yields an exact {\it weight complex} functor $t:\dmge\to K^b(\chowe)$;
 it could also be extended to  $t_{gm}: \dmgm\to K^b(\chow)$. In \S6
of \cite{mymot} it was also proved that these functors are conservative
(i.e. $t_{gm}(X)=0\implies X=0$).  

4.  $t$ gives $K_0(\dmge)\cong K_0(\chowe)$ and $K_0(\dmgm)\cong
K_0(\chow)$ (see \S6.4 of \cite{mymot}; a generalization and
certain variations of these results are described in \S\S5.3-5.5 of
\cite{bws}). See \S\ref{skz} below for the definitions of these groups and the formulation of the generalization.



5. Motivically functorial weight spectral sequences for any
cohomology theory $H:\dmge\to \au$ (generalizing Deligne's ones for
\'etale and  singular cohomology of varieties)
 were constructed (see \S6.6 and Remark 2.4.3 of \cite{bws};
 they were called
 Chow-weight spectral sequences since they correspond to the Chow weight structure; see \S\ref{swss}  below).

6. All triangulated subcategories and localizations of $\dmge$ were
 'described' (see \S8.1--8.2 of \cite{mymot}). In particular, one
obtains 'reasonable' descriptions of Tate motives and of the
(triangulated) category of birational motives (i.e. of the
localization of $\dmge$ by $\dmge(1)$; see \cite{kabir}) this way.

7. A certain category $\gd$ (of {\it comotives}) that contains
'nice homotopy limits' of
Voevodsky's motives was constructed (see \S3.1 and \S5 of \cite{bger}, and also \S\ref{scomot} below).
In particular, it contains certain (co)motives for all function fields
over $k$.

Some of the properties of $\gd$ are dual (in a certain sense) to the corresponding properties of the 'usual large
motivic' categories. In particular, though we have a covariant embedding
$\dmge\to\gd$, it yields a family of cocompact cogenerators for $\gd$.
This is why we call the objects of $\gd$ comotives.

Comotives are crucial for the proof of  the following
results.

8. There exist (certain) motivically functorial coniveau spectral sequences for cohomology
of arbitrary motives
  (see \S4.2 of \cite{bger}; cf.  \S7.4 of \cite{bws}).

 Besides, for $H$ represented by a motivic complex (i.e. an object of $\dme$)
we prove that these spectral sequences can be described in terms of the
homotopy $t$-truncations of $H$. This vastly extends seminal results
of Bloch and Ogus (see \cite{blog}).

9. Let $k$ be countable.

Then the cohomology of any smooth semi-local scheme (over $k$) is a
direct summand of the cohomology of its generic point; the
cohomology  of function fields contain twisted cohomology of their
residue fields (for all geometric valuations) as direct summands.

\end{theo}

\begin{rema}\label{rmmot}

1. The term 'weight complex' was proposed by Gillet and Soul\'e
in \cite{gs}. Their functor was essentially the restriction of
$t$ to motives with compact support of varieties (see \S6.6 of
\cite{mymot}). Besisdes, in \cite{gu} a functor that is essentially
$t\circ \mg$ was defined.

 Any of these functors allows to compute $E_2$ of the corresponding (Chow)-weight spectral sequences (see
 assertion 5).
Hence for (rational) singular/\'etale cohomology of varieties
(and motives) it
computes the factors  of the ('ordinary') weight filtration; whence the name. 

2. Parts 3--5 of the Theorem will be vastly generalized below
(to triangulated categories endowed with weight structures).

They follow from the existence of a certain {\it Chow weight
structure} for $\dmge$; whereas assertions 8--9 follow from the
existence of a certain {\it  Gersten weight structure} for some
triangulated $\gds$ such that $\dmge\subset\gds\subset\gd$.

3. Recently (independently in \cite{hebpo} and in \cite{brelmot}) it was also proved that the Chow weight structure could be defined for the category of Voevodsky's motives with rational coefficients over any 'reasonable' base scheme $S$ (in \cite{degcis} where the basic properties of
$S$-motives were established, they were called Beilinson's motives; one could either consider the 'large' category $DM(S)$ of $S$-motives or its subcategory $\dmcs$ of {\it constructible} i.e. 'geometric' objects here). The heart of this weight structure is 'generated' by (certain) motives of regular schemes that are projective over $S$ (tensored by $\q(n)[2n]$ for all $n\in \z$; see \S\ref{srmot} below). So, we obtain certain
analogues of parts 3--5 of the Theorem for $S$-motives also. 

In \S3 of \cite{brelmot} the weights for $S$-motives were also related with the 'classical' weights of mixed complexes of sheaves. To this end  the notion  of a {\it relative weight structure} was introduced; see Definition \ref{drwstr}
 below.

Besides, in \cite{lesm} for $S$ being a variety over a characteristic zero field an explicit differential graded description of a certain subcategory  of  $\dmcs$ 
was given; this is a generalization of assertion 1 of the Theorem. Possibly, this result could be extended to the whole $\dmcs$ (at least, with 
rational coefficients).

4. No explicit comparison functor
in the 'description' of part 1 is known (the two triangulated
categories in question are compared
by means of a third triangulated category).
Note also that the category of twisted complexes considered is a
'twisted' analogue
of $K^b(B)$ i.e. one considers morphisms and objects up to (a certain)
 homotopy equivalence.
 Hence in order to work with Voevodsky's motives one needs
constructions that do not depend on the choice of representatives
in  these homotopy equivalence
 classes. Weight structures really help here!

5. All the assertions of the theorem remain valid if we replace motives
with integral coefficients by those with rational (or $\znz$-) ones; see Remark \ref{rcoeff}(2).

Besides, the requirement $\cha k=0$ is only needed to apply the
resolution of singularities (that is required to prove some of the
statements in \cite{1}, which are necessary to deduce our results).
For motives with rational coefficients (we denote them by
$\chowe\q\subset \dmge\q\subset \dmgm\q$) for most of purposes it  suffices to
apply de Jong's alterations. In particular, this allows to prove the
'rational' analogues of assertions 3--5 also for any perfect $k$ of 
characteristic $p$. 

Moreover, a recent resolution of singularities result of Gabber (see Theorem 1.3 of \cite{illgab}) allows also to prove the
analogues of assertions 3--5 with $\zop$-coefficients (over $k$). Note here: Gabber's theorem could be called '$\zll$-resolution of singularities' (for all $l\in \p\setminus\{p\}$); yet weight structure methods allow to deduce motivic results with $\zop$-coefficients (that is a priori more difficult); see \cite{bzp}.

6. In \S6.3 of \cite{mymot} a certain length of motives was defined
(this is the 'length' of  $t(X)$). This is a motivic analogue
of the length of the weight filtration for mixed Hodge structures
(coming from cohomology of varieties). In particular, the length of
a motif
 of a smooth variety is is not greater  than its dimension and not
less than the length of the weight filtration for
its cohomology.

7. One can prove more than conservativity for $t$. In particular,
 $X\in \obj \dmgm$ is mixed Tate whenever $t_{gm}(X)$ is (see
 Corollary 8.2.3 of
\cite{mymot}).
\end{rema}

\section{Weight structures: basics}\label{sbws}

Now we define weight structures. They  are related with stupid
truncations of complexes (i.e. of objects of $K(B)$) in a way
similar  to the relation of $t$-structures  with canonical
truncations (see \cite{BBD} for the foundations of the theory of
$t$-structures); certainly, the distinctions here are also very
significant!

Stupid truncations are not very popular since they are not
canonical (whereas canonical truncations are canonical and functorial). Yet we will explain (starting from \S\ref{str} below) how
they do yield plenty of functorial ('cohomological') information; these results
are  new even for  $\cu=K(B)$. There are a lot of examples
when non-canonical constructions yield important functorial
information: projective and injective resolution of objects and
complexes over abelian categories allow to define derived functors;
nice compactifications and smooth hyper-resolutions of varieties
allow to define weight spectral sequences for \'etale and singular
cohomology; skeletal filtration for topological spectra allow to
construct Atiyah-Hirzebruch spectral sequences for their cohomology.
All of these observations have very natural 'explanations' inside
the theory of weight structures!

Weight structures have (at least) two distinct incarnations
important for Voevodsky's motives (related to weight and coniveau
spectral sequences), and also one that is relevant for the stable
homotopy category (in topology). Yet first we  illustrate some
basics of the theory on a (more) simple (though quite interesting)
example.

For $\cu=K(B)$
we denote by
$\cu^{w\le 0}$  the class of complexes, homotopy
equivalent to those concentrated in non-positive degrees;
we denote by $\cu^{w\ge 0}$ the class complexes,
equivalent to those concentrated in degrees $\ge 0$.

 Then the classes of complexes described satisfy the following
 properties (we write them
 down in the form that reminds the axioms of $t$-structures; this is
 quite convenient).

 \begin{defi}[Axioms of weight structures]\label{dwstr}

(i) $\cu^{w\ge 0},\cu^{w\le 0}$ are additive and Karoubi-closed in $\cu$.

(ii) {\bf 'Semi-invariance' with respect to translations.}

$\cu^{w\ge 0}\subset \cu^{w\ge 0}[1]$, $\cu^{w\le 0}[1]\subset
\cu^{w\le 0}$.

(iii) {\bf Orthogonality}.

$\cu^{w\ge 0}\perp \cu^{w\le 0}[1]$. 

(iv) {\bf Weight decompositions}.

 For any $X\in\obj \cu$ there
exists a distinguished triangle
\begin{equation}\label{ewd}
B[-1]\to X\stackrel{a}{\to} A\stackrel{f}{\to} B
\end{equation} such that $A\in \cu^{w\le 0}, B\in \cu^{w\ge 0}$.

For any triangulated category $\cu$ we will say that the classes
$(\cu^{w\le 0}, \cu^{w\ge 0})$ yield a weight structure if they
satisfy the properties listed.

\end{defi}

\begin{rema}\label{rbws1}
1. For $\cu=K(B)$ we can take weight decompositions coming 'stupid truncations' of complexes; see the illustration:
$$\begin{CD}
X=\dots @>{}>> X^{-2}@>{}>> X^{-1}@>{}>> X^{0}@>{}>> X^{1}@>{}>> X^{2}@>{}>>\dots
\\@.
 @.  @. @VV{a}V  @. @. @. \\
A=\dots @>{}>> X^{-2}@>{}>> X^{-1}@>{}>> X^{0}@>{}>> 0 @>{}>> 0 @>{}>> \dots
\\@.
 @.  @. @VV{f}V  @. @. @.
\\
B=\dots @>{}>> 0 @>{}>> 0 @>{}>> X^{1}@>{}>> X^{2}@>{}>> X^3@>{}>>\dots
\end{CD}$$

2. In this partial case ($\cu=K(B)$) we also have an opposite orthogonality
property ($\cu^{w\le 0}\perp \cu^{w\ge 1}$);  yet this additional orthogonality is not important, and
does not generalize to other (more interesting) examples.

3. For $t$-structures the orthogonality axiom is opposite; also, the arrows in
$t$-decompositions 'go in the converse direction'. These distinctions result in a drastic difference between the properties of these two types of structures. Note that dualization
does not change anything here (since the axiomatics of $t$-structures is
self-dual, as well as the one of weight structures).

4. We demand (in (i)) $\cu^{w\ge 0}$ and $\cu^{w\le 0}$ to be Karoubi-closed; this is  a technical condition that is not really important. The corresponding condition for $t$-structures is also true (though in contrast to the weight structure situation, it follows from the remaining axioms). 
\end{rema}

We also define the {\it heart} $\hw$ of $w$ (similarly to  hearts of
$t$-structures):
 $\obj\hw=\cu^{w=0}=\cu^{w\ge 0}\cap \cu^{w\le 0}$,
$\hw(X,Y)=\cu(X,Y)$ for $X,Y\in \cu^{w=0}$.

Now we list some very basic properties of weight structures (and
their hearts).

\begin{theo}\label{tbw}

\begin{enumerate}


\item\label{iext} $\cu^{w\le 0}$, $\cu^{w\ge 0}$, and $\cu^{w=0}$
are {\it extension-stable} i.e. for a distinguished triangle $A\to
B\to C$ if $A,C$ belong to $\cu^{w\le 0}$ (resp. to $\cu^{w\ge 0}$,
resp. to $\cu^{w=0}$) then $B$ belongs to the corresponding class
also.

    \item \label{isum}
    If $A\to B\to C\to A[1]$ is a distinguished triangle and
$A,C\in \cu^{w= 0}$, then $B\cong A\oplus C$.

\item $\hw$ is {\it negative}
i.e.  $\hw\perp \cup_{i>0}\hw[i]$. 

 \item \label{igen}
Conversely, for a triangulated $\cu$ let  an additive $D\subset \obj
\cu$ be negative; suppose that the smallest triangulated
 subcategory of $\cu$ containing $D$ is $\cu$ itself.
Then there exists a unique weight structure $w$ for $\cu$ such that
$D\subset \cu^{w=0}$; for it we have  $\hw=\kar_{\cu}D$ (see Theorem
4.3.2 of \cite{bws}).

One can construct all {\it bounded} weight structures (i.e. those ones
that  satisfy $\cap_{i\in \z} \cu^{w\le 0}[i] =\cap_{i\in \z}
\cu^{w\ge 0}[i]=\ns$) this way.
\end{enumerate}
\end{theo}

\begin{rema}\label{rbws2}
1. {\bf Examples}

Assertion \ref{igen} allows to construct the 'stupid' weight structure for
$K^b(B)$ mentioned above (note: as for $t$-structures, a single
$\cu$ may support more than one distinct weight structures).

Besides, in the stable homotopy category $SH$ there are no morphisms of
positive degrees between coproducts of the sphere spectrum $S^0$.
Hence assertion  \ref{igen} allows to construct a certain weight
structure for the subcategory $SH_{fin}\subset SH$ of finite spectra. In \S4 of
\cite{bws} several other existence of weight structures results (for
unbounded weight structures) were proved. In particular, they allow
to construct a certain $w_{S^0}$ for the whole $SH$ (see \S4.6 of
ibid.). The corresponding weight decompositions correspond to
cellular filtration of spectra; one can obtain Atiyah-Hirzebruch spectral
sequences this way (as {\it weight spectral sequences}; see below)!

 Lastly, $\chowe$ is negative inside $\dmge\subset \dme$; $\chow$ is
 negative
 inside $\dmgm$ (see 
 (\ref{evo})). This allows to construct certain
 {\it Chow}
 weight structures for all of these categories. We denote all of them
 by $w_{\chow}$, since they are compatible; see \S\S6.5-6.6 of
 \cite{bws}, and also 
 Remark
 \ref{rmmot} above.
 
 2. Assertion \ref{igen}  demonstrates that in the bounded case a weight structure could be completely described in terms of its heart; so instead of weight structures in this case one can consider only negative Karoubi-closed generating subcategories of $\cu$. Yet weight decompositions are very important (so it does not seem wise to avoid mentioning them in the axioms)!
 

 3. The obvious analogue of  assertion \ref{igen} for
 $t$-structures (i.e. we want to construct a $t$-structure such that a given
 positive $D\subset \cu$ lies in its heart) is very far
 from being true. So, negative subcategories of triangulated categories
 are much more valuable than positive ones! Besides, weight
 structures 'are more likely to exist for small triangulated categories'
 (than $t$-structures); see Remark 4.3.4 of \cite{bws}.

 4. Yet another distinction of weight structures from $t$-structures is
 demonstrated by assertion 3: distinguished triangles in $\cu$ do
 not yield non-trivial extensions in $\hw$.
 
 In fact, one may say that the notion of the heart of a weight structure is a 'triangulated analogue' of the category of projective (or injective) objects of an abelian category $\au$. Note here: we have $D(\au)(P,Q[i])=\ns$ if $i\neq 0$ and $P,Q$ are both projective (or injective) objects of $\au$; this allows to construct resolutions of objects of $\au$ (and hyperresolutions of complexes) that are functorial up to homotopy equivalence. The theory of weight structures demonstrates that one mostly needs $D(\au)(P,Q[i])=\ns$ if $i> 0$; the absence of 'positive extensions' is sufficient to prove certain functoriality of the corresponding 'resolutions' (i.e. Postnikov towers); see below. So, weight structures yield a vast generalization of projective and injective hyperresolutions!

 \end{rema}


 \section{On functoriality of weight decompositions;  truncations for
 cohomology}\label{str}

 Now we discuss to what extent weight decompositions are functorial,
 and how this allows to define nice canonical 'truncations' and
 filtration for cohomology.

 Weight decompositions (as in (\ref{ewd})) are (almost) never unique.
Still we will  denote any pair of $(A,B)$ as in (\ref{ewd}) by
$X^{w\le 0}$ and $X^{w\ge 1}$. $X^{w\le l}$ (resp. $X^{w\ge l}$)
will denote $(X[l])^{w\le 0}$ (resp. $(X[l-1])^{w\ge 1}$).
 $w_{\le i}X$ (resp. $w_{\ge i}X$) will denote $X^{w\le i}[-i]$
(resp.  $X^{w\ge i}[-i]$).

Now we observe that weight decompositions  are 'weakly functorial'.

\begin{pr}\label{pfwd}

1. Any $g\in\cu(X,Y)$ could be completed (non-uniquely) to a
morphism weight decompositions.

2. Moreover, for any $i\in\z$, $j>0$, $g$
extends
to a diagram
 \begin{equation}\label{efwd}
 \begin{CD} w_{\ge i+1} X@>{}>>
X@>{}>> w_{\le i}X\\
@VV{}V@VV{g}V@ VV{}V \\
w_{\ge i+j+1} Y@>{}>>
Y@>{}>> w_{\le i+j}Y \end{CD}
\end{equation}
in a unique way if we fix the corresponding weight decompositions.
\end{pr}

\begin{rema}

1. A nice illustration for assertion 1 is: for $\cu=\dmge$,
$w=w_{\chow}$, it  implies (in particular) that any morphism of
smooth varieties (coming from $\sv$, $\smc$, or $\dmge$) could be
completed in $\dmge$ to a morphism of (any choices of) their smooth
compactifications. Note: though one can prove this statement easily
without weight structures, yet it  is somewhat 'counterintuitive '.

2. For $\cu=K(B)$ assertion 2 means: if we fix the choice of weight
decompositions,  then the diagram
$$\begin{CD}
\dots @>{}>> X^{-2}@>{}>> X^{-1}@>{}>> X^{0}@>{}>> X^{1}@>{}>> X^{2}@>{}>>\dots
\\ @. @VV{g^{-2}}V @VV{g^{-1}}V @VV{g^{0}}V @VV{g^{1}}V @VV{g^{2}}V @. \\
\dots @>{}>> Y^{-2}@>{}>> Y^{-1}@>{}>> Y^{0}@>{}>> Y^{1}@>{}>> Y^{2}@>{}>>\dots
\end{CD}$$
is compatible with a unique  choice of the following diagram
$$\begin{CD}
(\dots @>{}>> X^{-2}@>{}>> X^{-1}@>{}>> X^{0})@>{f}>> (X^{1}@>{}>> X^{2}@>{}>>\dots)
\\ @. @VV{g^{-2}}V @VV{g^{-1}}V @. @VV{g^{1}}V @VV{g^{2}}V @. \\
(\dots @>{}>> Y^{-2}@>{}>> Y^{-1})@>{f'}>>( Y^{0}@>{}>> Y^{1}@>{}>> Y^{2}@>{}>>\dots)
\end{CD}$$
in $\cu$ (i.e. if we consider all morphisms up to homotopy equivalence).
\end{rema}

Proposition \ref{pfwd} immediately allows to construct some functorial
filtration and 'truncations' for cohomology (i.e. for some
contravariant $H:\cu\to\au$, that will usually be cohomological).

\begin{pr}
1. For any contravariant $H:\cu^{op}\to \au$, $j>0$, 
Proposition \ref{pfwd}(1) yields that the {\it weight filtration}
$W^iH(X)=\imm(H(w_{\le i}X)\to H(X))$ of $H(X)$ is
 $\cu$-functorial in $X$.

 2. Applying both parts of the proposition
 we obtain that
$H_1^i:X\mapsto \imm(H(w_{\le i}X)\to H(w_{\le i+j}X))$ also defines
a functor.

3. If $H$ is cohomological, $j=1$, $H_1^i$ is cohomological also.

4.  $H_2^i=\imm(H(w_{\ge i}X)\to H(w_{\ge i+1}X))$ is also
functorial  and cohomological (if $H$ is); there is a long exact
sequence of functors (i.e. it becomes a long exact sequence in $\au$
when applied to any object of $\cu$)
$$\dots \to H_2^i\circ [1]\to H_1^i\to H\to H_2^i\to H_1^i\circ [-1]\to\dots $$


\end{pr}

We call $H_1^i$ and $H_2^i$ {\it virtual t-truncations} of $H$. The
reason for this is that they  'behave as' if $H$ is 'represented' by
an object of some triangulated category $\du$, and the truncations
are 'represented' by its actual $t$-truncations with respect to some
$t$-structure of $\du$. We will observe that it is often the case in
the next section; yet note that virtual $t$-truncations can be
defined (and have nice properties) without specifying any $\du$ and
any $t$-structure for it (in fact, it is far from being obvious that
such $\du$ and $t$  exist always; even if they do, $\du$ is
definitely not determined by $\cu$ in a functorial way)!

Virtual $t$-truncations are studied (in detail) in \S2.5 of
\cite{bws} (there this concept was also developed for covariant functors; certainly, the difference is quite formal) and in \S\S2.3--2.5 of \cite{bger}. Also, ${\tilde W}_tH_n^{BM}(-,-)$  in Definition 5.8 of \cite{friha} are essentially (restrictions to motives of varieties of) virtual $t$-truncations of Borel-Moore homology with respect to $w_{\chow}$.

\section{Dualities of triangulated categories; orthogonal and
adjacent weight and $t$-structures}

Let $\du$  also be a triangulated category.

\begin{defi}\label{ddual}
1. We will call a (covariant) bi-functor
  $\Phi:\cu^{op}\times\du\to
\au$ a {\it duality} if  it is bi-additive, homological with respect
to both arguments; and is equipped with a (bi)natural transformation
$\Phi(X,Y)\cong \Phi (X[1],Y[1])$.

2. Suppose now that $\cu$ is endowed with a weight structure $w$,
 $\du$ is endowed with a $t$-structure $t$. Then we will say that $w$
 is (left) {\it orthogonal} to $t$ with respect to $\Phi$
 if the following
 {\it orthogonality condition} is fulfilled:
 \begin{equation}\label{edort}
 \Phi (X,Y)=0\text{ if: }X\in \cu^{w\le 0}
\text{ and }Y\in \du^{t \ge 1},\text{ or }X\in \cu^{w\ge 0}
\text{ and }Y\in \du^{t \le -1}.
 \end{equation} \end{defi}

\begin{rema}\label{rortog}
1. If $t$ is orthogonal to $w$, then: for any $X\in \cu^{w=0}$
the functor $Y\mapsto \Phi(X,Y)$ is exact when restricted to $\hrt$.

 Virtual $t$-truncations of $\Phi(-,Y)$ are 'represented' by
 $t$-truncations of $Y$: for example, $\Phi(X,Y^{t\ge i}[j])\cong \imm(\Phi([X^{w\ge -j},Y[i])\to
\Phi(X^{w\ge -1-j},Y[i-1]))$.

2. {\bf Adjacent structures}

A very important example of a duality is: $\du=\cu$, $\Phi(X,Y)=
\cu(X,Y)$. This duality is also {\it nice} (see Definition 2.5.1 of \cite{bger}); niceness is a technical
condition  
needed for  spectral sequences calculations (see below).

In this situation, we call orthogonal $w$ and $t$ {\it  adjacent
structures}; $w$ is (left) adjacent to $t$ whenever $\cu^{w\le
0}=\cu^{t\le 0}$; see \S4.4 of \cite{bws}.

3. {\bf Weight-exact functors; relation with adjoint functors}.

Recall now: if an exact functor $\cu\to \cu'$ is $t$-exact with
respect to some $t$-structures on these categories, its
(left or right) adjoint is usually not $t$-exact (it is only left
or right $t$-exact, respectively). This problem could be fixed 
if there exist
 adjacent weight structures for these $t$-structures (see Proposition
 4.4.5 of ibid.).

 We assume  that $\cu$ is endowed with a weight structure $w$ and its 
 adjacent $t$-structure $t$; $\cu'$ is endowed with a weight structure
 $w$ and its  adjacent $t$-structure $t$; $F:\cu\to \cu'$ is exact,
 $G:\cu'\to\cu$ is its left adjoint.

 We will say that $G$ is left (resp. right) weight-exact if
 $G(\cu'^{w\le 0})\subset \cu^{w\le 0}$ (resp.
 $G(\cu'^{w\ge 0})\subset \cu^{w\ge 0}$).

 Then: $G$ is left (resp. right) weight-exact whenever $F$ is right
 (resp. left) $t$-exact (in the well-known and similarly defined sense).

4. {\bf Examples}.

A simple example of adjacent structures is: if
$\operatorname{Proj}\au\subset \au$ denotes the full subcategory of
projective objects, $D^?(\au)$ (i.e. some version of $D(\au)$) is
isomorphic to the corresponding $ K^?(\operatorname{Proj}\au)$, then
for $\cu=D^?(\au)$  the canonic $t$-structure for $\cu$ is
adjacent to  the 'stupid' weight structure for $\cu\cong
K^?(\operatorname{Proj}\au)$ (mentioned above). Note that this
example allows to compute extension functors for $\au$ (and also
hyperextension ones i.e. morphisms in $D^?(\au)$)! Besides, the
spherical weight structure ($w_{S^0}$ for $SH$ mentioned above) is
adjacent to the Postnikov $t$-structure $t_{Post}$ (for $SH$).

Moreover, a process similar to the construction of
Eilenberg-Maclane spectra allows to construct a {\it Chow
$t$-structure} for $\dme$ such that  $\hrt_{\chow}\cong
\adfu(\chowe,\ab)$  (see \S7.1 of \cite{bws}). $t_{\chow}$ is
adjacent to the Chow weight structure for $\dme$; it is related with unramified cohomology (see \S7.6 of ibid.). Other related
calculations of hearts of orthogonal structures were made in
\S\S4.4--4.6 of ibid. and in \S6.2 of \cite{bger}.

Lastly, there also exists a nice duality $\gd^{op}\times \dme \to
\ab$ (see \S4.5 of \cite{bger}). If (the base field) $k$ is
countable, there also exists a triangulated category $\gds$ (such
that $\dmge\subset \gds\subset\gd$)  endowed with a {\it Gersten}
weight structure (see \S4.1 of ibid.),  that is orthogonal to the
homotopy $t$-structure for $\dme$ (defined in \S3 of \cite{1}). So, the
objects of its heart induce exact covariant functors from  $\hrt$
(i.e. the category of homotopy invariant sheaves with transfers) to
$\ab$. It is no surprise that this heart is 'generated' by comotives
of (spectra of) function fields (over $k$).

Note that in this case $\cu\neq \du$.

5. The recently proved Beilinson-Lichtenbaum conjecture
implies that the homotopy $t$-truncations of complexes of sheaves
that represent $\znz$-\'etale cohomology yield $\znz$-motivic
cohomology. Therefore one can express torsion motivic cohomology (of
smooth varieties, motives, and comotives) in terms of virtual
$t$-truncations of torsion \'etale cohomology with respect to the
Gersten weight structure. This allows to obtain some new formulae
for motivic cohomology; cf. \S\S7.4--7.5 of \cite{bws}  and Remark
4.5.2 of \cite{bger}.

 \end{rema}

\section{ Weight spectral sequences}\label{swss}

Applying $H$ to (shifted) weight decompositions of $X$ one obtains
an  exact couple $C_w(H,X)$ with:
$D_1^{pq}= H(X^{w\le -p}[-q])$,
$E_1^{pq}= H(X^{-p}[-q])$.

Here $X^i\in \cu^{w=0}$ are the terms of the {\it weight complex} of
$X$;  the latter coincides with $X$ for $\cu=K(B)$, was mentioned in
 Theorem \ref{main}(3) for $\cu=\dmge$ or $=\dmgm$, and will
be considered in \S\ref{swc} in the general case. We will call the
 spectral sequence corresponding to $C_w(H,X)$ a {\it weight spectral sequence} and
denote it by $T_w(H,X)$ (we will often omit $w$ in this notation).
 Under certain (quite weak) boundedness conditions this spectral sequence converges to $E_{\infty}^{p+q}T(H,X)=H(X[-p-q])$.
 Note that is natural to
denote $H(X[-i])$ by $H^i(X)$;
see also \S2.3--2.4 of \cite{bws} for more details.

$C_w(H,X)$ (and so also $T_w(H,X)$) is functorial
in $H$ (in the obvious way). Yet $C_w(H,X)$ (as well as $E_1(T_w(H,X))$) is not canonically determined by $X$ and $H$ (though any
$g\in \cu(X,X')$ could be extended to a morphism $C_w(H,X')\to C_w(H,X)$ for any choices of those). Still, such an extension is (almost) never  unique.

Yet this problem vanishes completely if one passes to the derived
exact couple! It is easily seen that $D_2$-terms are virtual
$t$-truncations of $H$ (defined in \S\ref{str} above); $E_2$ are
certain 'truncations from both sides'; so both are given by
cohomological functors $\cu\to\au$ (see loc.cit. and \S2.4 of
\cite{bger}). Hence $T(H,X)$ is (also) $\cu$-functorial (in $X$) starting from $E_2$.

Besides, the relation between virtual $t$-truncations and
truncations  with respect to an orthogonal $t$-structure (described
above) yields: for a nice duality $\Phi$, $H=\Phi(-,Y)$, $Y\in \obj
\du$, one has a functorial description of $T(H,-)$  (starting from
$E_2$)  in terms of $t$-truncations of $Y$; see Theorem 2.6.1 of
\cite{bger}. This is a powerful tool for comparing spectral
sequences (in this situation); it does not require constructing  any
complexes (and filtrations for them) in contrast to the method of
\cite{paran} (probably, originating from Deligne).

\begin{rema}[Examples; change of weight structures]\label{rwss}

1. Weight spectral sequences  generalize Deligne's weight spectral
sequences, coniveau, and Atiyah-Hirzebruch spectral sequences.

Weight spectral sequences corresponding to $w_{\chow}$ (we call them
{\it Chow-weight spectral  sequences} since they relate cohomology
of Voevodsky's motives with those of Chow motives) essentially
generalize Deligne's weight spectral sequences; see Remark 2.4.3 and
\S6 of \cite{bws}. For $H$ being \'etale or singular cohomology (of
motives) this yields motivic functoriality of $T_{w_{\chow}}(H,-)$
for integral (or torsion) coefficients. Note that the 'classical'
way of proving uniqueness of these spectral sequences uses Deligne's
weights for sheaves, and so requires rational coefficients (one also
uses heavily the fact that in this particular case weight  spectral
sequences degenerate at $E_2$).

One could also take the motivic cohomology theory for $H$. This yields
completely  new spectral sequences (yet see Remark 2.4.3(2) of
ibid.). This $T_{w_{\chow}}(H,-)$  does not degenerate at any
fixed level (even with rational coefficients, in general), and so
its functoriality definitely cannot be proved by 'classical'
methods.

2. Let $F:\cu\to \cu'$ be an exact functor that is right
weight-exact with respect to $w$ for $\cu$ and $w'$ for $\cu'$  (see
 Remark \ref{rortog}(3)); let $H:\cu'\to \au$ be cohomological. Then
in \S2.7 of \cite{bger} it was proved: for any $X\in \obj \cu$ there
exists some comparison morphism of
 weight spectral sequences $M:T_{w}(H\circ F,X) \to T_{w'}(H,F(X))$.
 Moreover, this morphism is unique and additively functorial
  starting from $E_2$.
The proof uses a natural (and easy) generalization of (\ref{efwd}).

In particular, this yields  comparison functors from Chow-weight  spectral sequences  to coniveau ones (cf. \S\ref{sgbir} below for more details).

If $F$ is left weight-exact, there exists a comparison
transformation  $N$ in the inverse direction. We call both $M$ and
$N$ 'change of weight structures' transformations.

3. Using the Gersten weight structure (for $\gds$, see above) one
can extend 'classical' coniveau spectral sequences from (motives of) smooth varieties to $\gds\supset \dmge$ in a
natural way (for an arbitrary cohomology theory $H$ that factorizes through 
$\dmge$, such that $\au$ satisfies AB5). This also yields motivic
functoriality of coniveau spectral sequences (which is far from
being obvious from their definition; see Remark 4.4.2 of
\cite{bger}). Note also that we obtain this functoriality for a not
necessarily countable $k$, since one can always define the coniveau
spectral sequence for $(H,X)$ over $k$ as the limit of the related
coniveau spectral sequences over countable perfect fields of
definition of $X$ (see \S4.6 of ibid.). Here we use the 'change of
weight structure' transformations (that we denoted by $N$ above).

The orthogonality of the Gersten weight structure  with the homotopy
$t$-structure (for $\dme$; see the previous section) yields that the
coniveau spectral sequence for $H$ represented by some $Y\in \obj
\dme$ could be described in terms of the homotopy $t$-truncations of
$H$. This extends vastly the
 coniveau spectral sequence calculations of Bloch\ \&\ Ogus
(in \cite{blog}; see \S4.5 of \cite{bger}).

4. Since $t_{Post}$ and $w_{S^0}$ are adjacent, we obtain  the
well-known fact: the Atiyah-Hirzebruch spectral sequence converging to $[X,Y]$ for
$X,Y\in \obj SH$ could be expressed either in terms of the
$t_{Post}$-truncations of $Y$ or in terms of $w_{S^0}$-truncations
of $X$ (i.e., in terms of cellular filtration of $X$).

\end{rema}

\section{More on weight structures} \label{sore}

\subsection{'Functoriality' of weight structures: localizations  and gluing}

Weight structures could be carried over to localizations and
also 'glued' similarly to $t$-structures. 

If $w$ (for $\cu$) induces a weight structure also on some
triangulated $\du\subset \cu$, then it also induces a compatible
weight structure on the Verdier quotient $\cu/\du$; its heart could
be easily described (in terms of the hearts of $\cu$ and $\du$ in a way that is quite distinct from those for $t$-structures; see \S8.1 of \cite{bws}). 

Moreover, one can
glue weight structures (i.e. recover a weight structure for $\cu$
from  those for $\du$ and $\cu/\du$ when certain adjoint functors
exist) in a way that is just slightly different from those for
$t$-structures (see \S8.2 of ibid.). We discuss an
interesting example of such a gluing in  \S\ref{sgbir} below. 

 This statement was also used in \S2.3 of \cite{brelmot} in  (one of the methods  of) the construction
of the Chow weight structure  for  motives over $S$. 

\subsection{The weight complex functor}\label{swc}

There are two ways to construct the weight complex functor for a
general $(\cu,w)$ (that generalizes the exact conservative functor
$t:\dmge\to K^b(\chowe)$ mentioned in Theorem \ref{main}).

First we describe the 'rigid' method. Suppose that $\cu$ has a
'description' in terms of {\it twisted complexes} over a negative
differential graded category (i.e. a {\it differential graded
enhancement}; see \S2 of \cite{mymot} or \S6 of \cite{bws}). Suppose
also that $w$ is compatible with  this enhancement (i.e. that $w$
coincides with the weight structure given by Proposition 6.2.1 of
ibid.). Then there exists an exact weight complex functor $t:\cu\to
K(\hw)$; see \S6.3 of ibid. (actually, in loc.cit. only bounded
twisted complexes are considered, so the target
of $t$ is $K^b(\hw)$). 

The main disadvantage of this method is that it requires some  extra
information on $\cu$. A differential graded enhancement does not
have to exist at all (for a general $\cu$; for example,  $SH$ has no differential graded enhancements); an exact functor does
not have to extend to enhancements (and if such an extension
exists, it is not necessarily unique).

Luckily, in \cite{bws}  another method was developed; it always
works and does not depend on any extra structures.
 There is a construction that associates a certain complex to
 each $X\in \obj\cu$ for any $\cu$ and depends only on $w$. It
 is closely related with the definition of a {\it weight Postnikov tower}
 for $X$ (see Definitions 1.1.5 and 2.1.2 of \cite{bger}). The terms of
 the (weight) complex $t(X)$ are
 $X^i=\co (w_{\le i-1}X\to w_{\le i}X)[i] \cong \co (w_{\ge i}X
 \to w_{\ge i+1}X)[i-1]$ (see Remark 2.1.3 of loc.cit.);
 the corresponding triangles yield some boundary morphisms
 $X^i\to X^{i+1}$ (see \S2.2 of \cite{bger}). It is easily seen that any
 $g\in \cu(X,X')$ is compatible with some  $t(g):t(X)\to t(X')$. This
 method has the following serious disadvantage: 
in general, $t(g)$ is only well-defined up to morphisms of the form
$df+gd$ (i.e. modulo an equivalence relation that is more coarse
than homotopy equivalence of morphisms of complexes). Still this
equivalence relation has certain nice properties: equivalent morphisms yield the same map on 
the cohomology of complexes;
the homotopy
equivalence class of $t(X)$ does not depend on the choices mentioned.
So, we obtain a certain {\it weakly exact} functor $\cu\to
K_\w(\hw)$ (see Definition 3.1.5 of loc.cit.). For any $H$ one has
$E^{pq}_1T(H,X)=H(X^{-p}[-q])$; hence $E^{**}_2T(H,X)$ can be
described in terms of $t(X)$
(in a functorial way); see Remark 3.1.7 of loc.cit. 

In the case $\cu=SH$ we have $K_\w(\hw)=K(\hw)$; so $t$ is actually
an exact functor  (see Remark 3.3.4 of ibid.).

Moreover, this ('weak') weight complex functor is compatible with
the 'strong' one given by the differential graded approach; see
\S6.3 of ibid.
It  is conservative if $w$ is {\it bounded}  (i.e. if
$\cap_{i\in \z} \cu^{w\le 0}[i] =\cap_{i\in \z} \cu^{w\ge
0}[i]=\ns$); see Theorem 3.3.1 of ibid. for the proof of this fact
and of several other nice properties of $t$.

\subsection{Certain $K_0$-calculations}\label{skz}

Suppose that $w$ is bounded, $\hw$ is idempotent complete.
Then  $\cu$ is idempotent complete also; see Lemma 5.2.1 of ibid. In
particular, this yields that $\dmge$ is generated by
$\chowe$ (i.e. the only  strict full triangulated subcategory of
$\dmge$ containing $\chowe$ is $\dmge$ itself); it seems that \S3.5
of \cite{1} does not contain a complete proof of this statement.

Besides, we have $K_0(\cu)\cong K_0(\hw)$. Recall that the generators of $K_0(\cu)
$ (resp. $K_0(\hw)$) are $[X]$, $X\in \obj \cu$ ($X\in \obj \hw$),
and the relations are: $[B]=[A]+[C]$ if $A\to B\to C$ is a
distinguished triangle (resp. $B\cong A\bigoplus C$).

In particular, we obtain Theorem \ref{main}(4) this way. 

\subsection{A generalization: relative weight structures}

Now we describe a formalism that generalizes those of weight structures. It is actual since in the (derived) category of mixed complexes of sheaves 
over a variety $X_0$ defined over a finite field $\fq$ the subcategories of objects of non-positive and non-negative weights do not quite satisfy the orthogonality axiom (iii) of Definition \ref{dwstr}. So, we adjust this axiom in order make it compatible with Proposition 5.1.15 of \cite{BBD}. Note here: in our notation the roles of $\cu^{w\le 0}$ and $\cu^{w\ge 0}$ are permuted with respect to the notation of (\S5.1.8 of) ibid.

\begin{defi}\label{drwstr}

Let $F:\cu\to\du$ be an exact 
 functor (of triangulated categories).

A pair of extension-stable (see Theorem \ref{tbw}(\ref{iext})) Karoubi-closed 
subclasses $\cu^{w\le 0},\cu^{w\ge 0}\subset\obj \cu$ for
a triangulated category $\cu$ will be said to define a relative weight
structure $w$ for $\cu$ with respect to $F$ (or just and $F$-weight structure)
if they  satisfy  conditions (ii) and (iv) of Definition \ref{dwstr}, as well as the following orthogonality assumptions:

$\cu^{w\ge 0}\perp \cu^{w\le 0}[2]$;
$F$ kills all morphisms between 
$\cu^{w\ge 0}$ and $\cu^{w\le 0}[1]$. 
\end{defi}

Relative weight structures satisfy several properties similar to those of 'absolute' weight structures (note: an 'absolute' weight structure is the same thing as an $\id_{\cu}$-weight structure); see below.

\section{'Motivic' weight structures; comotives;
gluing Chow and Gersten structures from 'birational slices'}\label{swmot}

We
briefly summarize how
 weight structures help in the proof of Theorem \ref{main}
 (this information could be found above, yet it is somewhat scattered).
 We also make several other remarks.

 As explained above, weight structures yield a mighty instrument for constructing and
 studying certain functorial spectral sequences for cohomology functors
 (defined on a triangulated category $\cu$); so they also yield certain
 functorial ('weight') filtration. They also describe how objects of
 $\cu$ could be 'constructed from' objects of a 'more simple'
additive $\hw\subset \cu$.

We have two main 'motivic' weight structures. They correspond to
(Chow)-weight and coniveau spectral sequences, respectively. Note
that both of these spectral sequences were 'classically' defined
only for cohomology of varieties; still our approach allows to
define them for arbitrary Voevodsky's motives, and also yields their
motivic functoriality (which is very far from being obvious).

\subsection{Chow weight structure(s); relation with the  motivic
$t$-structure and weight filtration}

 Our first ('motivic') weight structure (being more precise, we have
 a system of compatible weight structures on various 'motivic' categories) is $w_{\chow}$; it is defined on $\dmge\subset \dmgm$, its heart is $\chowe\subset \chow$; $w_{\chow}$ can also be extended to $\dme$ and $\gd$. So, it closely relates $\dmge$ with $\chowe$ (in particular, the weight complex functor $\dmge\to K^b(\chowe)$ is conservative; note that $\dmge$ is very far from being isomorphic to $K^b(\chowe)$!).  So, the cohomology of Voevodsky's motives can be 'functorially related' with the cohomology of Chow ones; one obtains a vast generalization of Deligne's weight spectral sequences.
 
 Besides, there exists a {\it Chow
$t$-structure} for $\dme$ such that  $\hrt_{\chow}\cong
\adfu(\chowe,\ab)$; $t_{\chow}$ is
adjacent to the Chow weight structure for $\dme$.

Now we relate $w_{\chow}$ with the 'usual expectations for weights of
motives';  see \S8.6 of \cite{bws} for more details.

Conjecturally, $\dmge\q$ (and $\dmgm\q$) should  support a certain ('mixed')
motivic $t$-structure ($t_{MM}$, whose heart is the abelian category
$MM^{eff}\subset MM$ of mixed motives) and a  {\it weight filtration} (by certain
triangulated subcategories); the latter one comes from  certain
weight filtration  functors $MM\to MM$ (compatible via cohomology
with the weight filtration of mixed Hodge structures and of mixed
Galois modules; these functors are idempotent). So, there should be
three important filtrations for $\dmge\q\subset\dmgm\q$ altogether.

Now, one can easily verify that the (widely believed to be true, yet conjectural)
properties of the two conjectural filtrations mentioned yield: for
a subcategory of objects that are 'pure of some fixed weight $i$' with respect to
any one of the three filtrations mentioned, the filtrations induced by two
remaining structures differ only by a shift of indices (that depends on $i$). In particular, $t_{MM}$
'should split' Chow motives into components that are 'pure with
respect to the weight filtration'; $w_{\chow}$-weight decompositions
induce the (conjectural!) weight filtration for mixed motives.
Note here: though weight decompositions (of objects of triangulated categories) are (usually) highly non-unique, for any $i\in
\z$, $X\in \obj MM\subset \obj \dmge\q$, there 'should exist' a unique
weight decomposition of $X[i]$ such that $w_{\le i}X,\, w_{\ge
i+1}X\in MM$; this choice of $w_{\ge i+1}X$ is what one expects to
be the corresponding level of the weight filtration of $X$ in $MM$.

In \cite{wildat} this (conjectural) picture  was justified in the
case when $k$ is a number field for the triangulated category
$DAT\subset\dmge\q$ (of so-called  Artin-Tate motives; this is the
triangulated subcategory of $\dmge\q$ generated by Tate twists of
motives of spectra of finite extensions of $k$). It was also
 shown that
the restriction of $w_{\chow}$ to $DAT$ can be completely characterized
 in terms of weights of singular homology. Actually, this corresponds
 to the fact that the triangulated category $DHS$ of mixed Hodge
complexes has a weight filtration (by triangulated subcategories)
and  could be endowed with a weight structure; these filtrations and
the 'canonical' $t$-structure for $DHS$ are connected by the same
relations as those that 'should connect' the corresponding
filtrations of $\dmge\q\subset \dmgm\q$. Besides, it could be easily seen
that singular (co)homology is weight-exact. 

\subsection{Comotives; the Gersten weight structure}\label{scomot}

 Our second 'motivic' weight structure is the Gersten weight
 structure $w$ defined on the category $\gds\supset \dmge$ (for a
 countable $k$).
Here $\gds$ is a full triangulated subcategory
of a certain category $\gd$ of {\it comotives} (already mentioned in
Theorem \ref{main}).

The idea is that $w$ should be orthogonal to the homotopy
$t$-structure  on $\dme$ (recall that the latter is the restriction
of the canonical $t$-structure of the derived category of Nisnevich
sheaves with transfers). So, $\hw$ is 'generated' by comotives of
function fields over $k$ (note that these are Nisnevich points).

It follows that  $w$ cannot be defined on $\dmge$ (or on $\dme$).
The problem with $\dme\supset \dmge$ is that there are no 'nice'
homotopy limits in it. In order to have them one needs 'nice'
(small) products; one also needs the objects of $\dmge$ to be
cocompact (in this 'category of homotopy limits'). $\dme$ definitely
does not satisfy these conditions. Instead in \S5 of \cite{bger} a
category $\gdp$ that is opposite to a certain category of {\it
differential graded modules} (i.e. covariant differential graded
functors from the {\it differential graded enhancement} of $\dmge$
to complexes of abelian groups) was considered; $\gd$ is its
homotopy category (with respect to a certain closed model structure;
so it is opposite to the corresponding derived category of
differential graded modules). So, we have a contravariant Yoneda
embedding of $\dmge$ to the category opposite to $\gd$ whose image
consists of compact objects; in this category 'nice' homotopy
colimits exist. Thus, inverting  arrows we obtain a 'nice' category
of comotives. Inside $\gd$ we define $\gds$ as  its smallest
Karoubi-closed triangulated category that contains (countable) products of comotives of
functions fields. Note: we need $k$ to be countable since without
this the author does not know how to prove that (our candidate for)
$\hw$ is negative; still comotives can be defined over any perfect
$k$.

The general theory of weight spectral sequences yields those for
cohomological  functors $\gds\to \au$. A (minor)  problem here is that
$\gds$ is 'large'; 
yet any
$H:\dmge\to \au$ has a 'nice' extension to $\gds$ (and also to
$\gd\supset \gds$) if $\au$ satisfies AB5 (see Proposition 4.3.1 of
\cite{bger}).  So, we can consider weight spectral sequences
$T=T_{w}(H,X)$ for any such $H$ and any $X\in \obj \dmge$ (or $X\in
\obj \gds$). It turns out that for $X$ being the motif of a smooth
variety,  $T$ is isomorphic to the coniveau spectral sequence
(corresponding to $H$) starting from $E_2$; see Proposition 4.4.1 of
ibid. So, we call $T$ a coniveau spectral sequence for any $X$. As
in the case of   'classical' coniveau spectral sequences, if $H$ is
represented by an object of $\dme$, $T_w(H,X)$ can be described in
terms of cohomology of $X$ with the coefficients in the homotopy
$t$-truncations of $H$ (see Corollary 4.5.3 of ibid.); this fact
extends the related results of Bloch-Ogus and Paranjape (see
\cite{blog} and \cite{paran}). Our latter result follows from the
existence of  a nice duality $\gd^{op}\times \dme \to \ab$.

\begin{rema}

$w$ can be restricted to the category $DAT\subset \dmge$ of
Artin-Tate motives (mentioned above; one may take integral
coefficients here; $k$ is any perfect field). Indeed, we don't need
comotives here, since (co)motives of (spectra of) finite extensions
of $k$ belong to $\obj \dmge$.

We explain this in more detail. $DAT$ is generated by
$\mg(F)(j)[j]$, where $F$ runs through all  (spectra of) finite
field extensions of $k$, $j\ge 0$.
$D=\{\oplus_i\mg(F_i)(j_i)[j_i]\}$ is a negative (additive)
subcategory of $DAT$, so Theorem \ref{tbw}(4) implies: there
exists a weight structure $w_{DAT}$ with $D\subset \hw_{DAT}$. Since
$\hw_{DAT}\subset \hw(\subset \gd)$, we obtain that $w_{DAT}$ is
compatible with $w$ (at least, for a countable $k$).

In particular, this implies that coniveau spectral sequences for
cohomology of any $X\in \obj DAT$ have quite 'economical'
descriptions (starting from $E_2$).
\end{rema}

\subsection{Comparison of weight structures; 'gluing from birational
slices'}\label{sgbir}

First we describe the relation between $T'=T_{w_{\chow}}(H,X)$ and
$T=T_{w}(H,X)$ (for $X\in \obj \dmge\subset \obj\gds$).
The 'change of weight structure transformation' (see Remark
\ref{rwss}) yields some morphism $M:T\to T'$ (functorially starting
from $E_2$; see \S4.8 of \cite{bger}). $M$ is an isomorphism if $H$
is {\it birational} i.e. kills $\dmge(1)$; here $-\otimes \z(1)$ is
the Tate twist isomorphism of $\dmge$ into itself.

Now,  $-\otimes \z(1)$ can be extended from $\dmge$ to $\gd$ (see
\S5.4.3 of ibid.); this is also true for $w_{\chow}$ (see \S4.7 of
ibid.). It is easily seen that $w$ and $w_{\chow}$ induce the same
weight structure $w_{bir}$ on the category of {\it birational
comotives} $\gdbr=\gd/\gd(1)$ (the Verdier quotient); the heart of
this localization contains images of all (co)motives of all smooth
varieties. One obtains that  (roughly!) $w$ and $w_{\chow}$
'coincide on slices' and only differ by the value of a single
integral parameter: $w$ is $-\otimes \z(1)[1]$-stable  and
$w_{\chow}$ is $-\otimes \z(1)[2]$-stable!

We try to make this more precise; see \S4.9 of ibid. for more
details.  We consider  the localizations $\gd/\gd(n)$ for all $n>0$.
Though none of them is isomorphic to $\gd$, they 'approximate it
pretty well'. 
Also, for any $n$ we have  a short exact sequence of
triangulated categories
  $\gd/\gd(n)\stackrel{i_*}{\to} \gd/\gd(n+1)\stackrel{j^*}{\to} \gdbr$.
   Here the notation for functors comes from the 'classical' gluing
   data setting (cf. \S8.2 of \cite{bws}); $i_*$
  can be given by $-\otimes\z(1)[s]$ for any $s\in\z$, $j^*$ is just
   the localization. Now, if we
  choose $s=2$ then
  both
  $i_*$ and $j^*$ are
 weight-exact with respect to weight structures induced by $w_{\chow}$
 on the corresponding categories;
 if we choose $s=1$ these functors are weight-exact with respect to
 the weight structures coming from $w$. So, the Chow and Gersten weight
 structures induce  weight
structures on the localizations $\gd(n)/\gd(n+1)\cong \gdbr$ (we call
these localizations 'slices') that differ only by a shift.

One can show that for any
  short exact sequence
  $\du\stackrel{i_*}{\to} \cu\stackrel{j^*}{\to} \eu$ of
  triangulated categories,
 if $\du$ and $\eu$ are endowed with weight structures, then there
  exist at most one weight structure on $\cu$ such that both $i_*$
  and $j^*$ are weight-exact. So, if one
    calls the filtration of $\gd$ by $\gd(n)$ the {\it slice filtration}
(this term was already used by A. Huber, B. Kahn, M. Levine,
V. Voevodsky, and other authors for other 'motivic categories'), then
one may say that the weight structures induced by $w$ and $w_{\chow}$
on all $\gd/\gd(n)$ 'can be recovered from slices'; the  only
difference between them is 'how we shift the slices'!

Moreover, Theorem 8.2.3 of \cite{bws} shows that if both adjoints to
both $i_*$
  and $j^*$ exist,
  then one can use this gluing data in order to 'glue' (any pair) of
   weight structures for $\du$ and $\eu$ into a weight structure for
   $\cu$. So, suppose that we have a weight structure $w_{n,s}$ for
$\gd/\gd(n)$ that is $-\otimes (1)[s]$-stable and 'compatible with
$w_{bir}$ on all slices'.  Then we can also construct $w_{n+1,s}$
satisfying similar properties, since general homological algebra
yields that all adjoints needed exist in our situation. So,
$w_{n,s}$ exist for all $n>0$ and all $s\in \z$. Moreover, there exists a 'large' subcategory of $\gd$ (containing $\dmge$) that for any $s$
can be endowed with a weight structure $w_s$ compatible with all $w_{n,s}$.
     Hence Gersten and
Chow weight structures (for $\gds/\gds(n)\subset \gd/\gd(n)$) are
members of a rather natural family of weight structures indexed by a
single integral parameter! It could be interesting to study other
members of this family (for example, the one that is $-\otimes
\z(1)$-stable).

\subsection{Weights for relative motives and mixed sheaves}\label{srmot} 

Let $S$ be a scheme of finite type over some excellent noetherian scheme $S_0$ of dimension $\le 2$. 

As we have already said, on the category $\dmcs$ (of {\it constructible} i.e. 'geometric' motives with rational coefficients over $S$) there exists a weight structure $\wchow$ whose heart $Chow(S)$ is the idempotent completion of $\{p_!\q_P(n)[2n]\}$, for $p:P\to S$ running through projective (or proper) morphisms such that $P$ is regular, $n\in \z$ (see \S3 of \cite{hebpo} and  \S2.1 of \cite{brelmot}). 

The corresponding Chow-weight spectral sequences yield: for any cohomological $H:\dmcs\to \au$, $X\in \obj \dmcs$, there exists a filtration on $H^*(X)$ (that is $\dmcs$-functorial in $X$) whose factors are subfactors of cohomology of some regular projective  $S$-schemes'; see Remark 3.3.2(3) of \cite{brelmot}.
Besides, the (Chow)-weight filtration of cohomology yields a natural way of description of the 'integral part' of the motivic cohomology of a variety over a number field (as constructed in \cite{scholl}; see Remark 3.3.2(4) of \cite{brelmot}).

We also obtain that $K_0(\dmcs)\cong K_0(Chow(S))$ (cf. \S\ref{skz}), and define a  certain 'motivic Euler characteristic' for  $S$-schemes (in \S3.2 of \cite{brelmot}). The author hopes that these results could be useful for motivic integration. 

Now denote by $\mathbb{H}$  the \'etale realization functor $\dmcs\to \dsh$, where $\dsh=\dsh(S)$ is the category $\dbm(S,\ql)$ of mixed complexes of $\ql$-\'etale sheaves
as considered in \cite{huper} and in \cite{BBD}.  
Then $\mathbb{H}$ sends Chow motives over $S$ to pure complexes of sheaves (see Definition 3.3 of \cite{huper} and \S3.4 of \cite{brelmot}). We deduce certain consequences from this fact.


Suppose that $S$ is a finite type $\spe \z$-scheme. 
We take $H_{per}$ being the perverse \'etale cohomology theory  i.e. $H_{per}^i(M)$ (for $M\in \obj \dmcs$, $i\in\z$) is the $i$-th cohomology of $\mathbb{H}(M)$ with respect to the perverse $t$-structure of $\dsh$ (see  Proposition 3.2 of \cite{huper}). 
 Then 
 $T_{\wchow}(H_{per},M)$ for any $M\in \obj \dmcs$ yields: all $H_{per}^i(M)$ have weight filtrations (defined using Definition 3.3 of loc.cit., for all $i\in\z$). Note that this is not at all automatic (for perverse sheaves over $S$); see Remark 6.8.4(i) of \cite{janbook}. Certainly, one can replace perverse sheaves over $S$ here by $\ql$-adic representations of the absolute Galois group of the function field of $S$; cf. \S6.8 of loc.cit.



Now let $S=X_0$ be a variety over a finite field $\fq$; let $X$ denote $X_0\times_{\spe \fq}\spe \ff$, where $\ff$ is the algebraic closure of $\fq$. 
The results of \S5 of \cite{BBD} (along with some of the results of \cite{brelmot}) yield that the category $\dsh(=\dbm(X_0,\ql))$  can be endowed with an $F$-weight  structure $w_{\dsh}$ whose heart is the category of pure complexes of sheaves,  for $F$ being the extension of scalars functor $\dsh\to D^b(X,\ql)$; see Proposition 3.6.1 of \cite{brelmot}. 
In particular, we obtain that any object $M$ of $\dsh$ possesses a 'filtration' (a {\it weight Postnikov tower}) whose 'factors' belong to $\hw_{\dsh}$. 

 Next, our $\mathbb{H}$ is a {\it weight-exact functor}
 (i.e. it sends $\dmcs^{w_{\chow}\le 0}$ to $\dsh^{w_{\dsh}\le 0}$ and sends $\dmcs^{w_{\chow}\ge 0}$ to $\dsh^{w_{\dsh}\ge 0}$). 
Hence this is no wonder that the weight-exactness properties of motivic base change functors (for $DM^c(-)$; see Proposition 3.8 of \cite{hebpo}, and  Theorem 2.2.1 and Proposition 2.3.4 of \cite{brelmot}) are parallel to the 'stabilities' 5.1.14 of \cite{BBD}.

Lastly, let  $G:D^b(X,\ql)\to \au$ be any cohomological functor, $H=G\circ F$, $M\in \obj \dmcs$. Then the weight-exactness of $\mathbb{H}$ yields that the (Chow)-weight filtration  for $(H\circ \mathbb{H})^*(M)$ is exactly the $w_{\dsh}$-weight filtration for $H^*(\mathbb{H}(M))$; cf. Proposition 3.5.5(II2) of \cite{brelmot}.

Very probably, some  analogues of these results are valid for  $\mathbb{H}$ replaced by a 'Hodge module' realization of motives (for $S$ being a complex variety); the problem is that (to the knowledge of the author) no such realization is constructed at the moment.

 \section {Possible applications to finite-dimensionality of motives }

Recall that $\dmge\subset \dmgm$, as well as their 'rational
versions'  $\dmge\q\subset \dmgm\q$ (see  Remark
\ref{rcoeff}(2)) are tensor triangulated categories. This allows to
define external and symmetric powers of objects in two latter
categories, since those are direct summands of tensor powers (for $\q$-linear motivic categories).

 $M\in \dmgm\q$  is called {\it Kimura-finite} (or finite-dimensional)
 if $M=M_1\bigoplus M_2$, where some external power of $M_1$ and some
 symmetric power of $M_2$ is $0$. In this case $M_1$ is called
 {\it evenly finite-dimensional}.
  Now, $t\q:\dmge\q\to K^b(\chowe\q)$ (the rational version of the
  weight complex) is a
conservative tensor functor; so   $X\in \obj \dmge$ (or $\dmgm$, or
 $\obj \dmgm\q$) is Kimura-finite whenever
$t_{gm}\q(X)$ is.

Now we describe a series of motives that 'should be'  finite-dimensional  (very similar
objects were considered by A. Beilinson and M. Nori though in somewhat
 distinct contexts).

Let $X/k$ be smooth affine of dimension $n$, $Y$ be its generic
 hyperplane section (with respect to some projective embedding).
Then for $M=(Y\to X)$  the only non-zero cohomology is
$H_{et}^n(M_{k^{alg}})$. Hence
some external power of $M\otimes \q[-n]$  'should' vanish
(since a certain external
 power of its cohomology vanishes). So
$M[-n]$ 'should be' evenly finite-dimensional. We can also pass to
$K^b(\chowe\q)$ here (i.e. consider $t(M)$ instead of $M$) since
the rational version of the weight complex functor is a tensor
functor.

\begin{rema} 1. If all such $M$ are  Kimura-finite at least
numerically (i.e. we consider their images in
$K^b(\operatorname{Mot_{num}})$ obtained via $t$),
 then  one can prove that $\operatorname{Mot_{num}}$ is a tannakian
 category.

2. Widely-believed conservativity of \'etale cohomology (as a
functor on $\dmge\q$) immediately implies that all such $M$ are
Kimura-finite indeed (as mentioned above). Alternatively, it is
possible to deduce Kimura-finiteness of $M$ from a certain  weak
Lefschetz for motivic cohomology. The latter 'should be true' since
it easily follows  from the (widely believed, yet conjectural!)
existence of a 'reasonable' motivic $t$-structure for $\dmge\q$.

Unfortunately,
 the author has no idea how to prove anything here unconditionally.
\end{rema}

\end{document}